\documentstyle[draft%
,amscd%
,syntonly%
,amssymb%
,psfig%
]{amsart}

\newtheorem{thm}{Theorem}   
\newtheorem{cor}{Corollary}
\newtheorem{lem}{Lemma}

\newtheorem{defn}{Definition}

\begin{document}



\title[large volume genus 2]
{Genus 2 closed hyperbolic 3-manifolds of arbitrarily large volume}
\author{Jennifer Schultens}
\address{1784 N Decatur Rd, Suite 100 \\
	Dept of Math and CS \\
	Atlanta, GA 30322}
  
\email{jcs@@mathcs.emory.edu}
\thanks{Research partially supported by NSF grant DMS-9803826}
\maketitle

\begin{abstract}
We describe a class of genus 2 closed hyperbolic 3-manifolds of
arbitrarily large volume.
\end{abstract}

\vspace{5 mm}

The purpose of this note is to advertise the existence of a class of
genus 2 closed hyperbolic 3-manifolds of arbitrarily large volume.
The class described here consists merely of appropriate Dehn fillings
on $2$-bridge knots.  That this class has the properties claimed
follows directly from \cite{CL}, \cite{L}, and the Gromov-Thurston
$2\pi$-Theorem.  The existence of such a class of hyperbolic
$3$-manifolds is known, as pointed out by Daryl Cooper \cite{C}, who
mentions that branched covers of the figure $8$ knot provide another
such class.  We believe that the existence of such a class deserves to
be more widely known.  For general definitions and properties
concerning knot theory, see \cite{Lic} or \cite{Liv}.

The following definition and theorem are due to M. Lackenby.

\begin{defn}
Given a link diagram $D$, we call a complimentary region having two
crossings in its boundary a \underline{bigon region}.  A \it{twist} is
a sequence $v_1, \dots, v_l$ of vertices such that $v_i$ and $v_{i+1}$
are the vertices of a common bigon region, and that is maximal in the
sense that it is not part of a longer such sequence.  A single
crossing adjacent to no bigon regions is also a twist.  The \it{twist
number t(D)} of a diagram $D$ its number of twists.
\end{defn}

\begin{thm}
(Lackenby) Let $D$ be a prime alternating diagram of a hyperbolic link
$K$ in $S^3$.  Then $v_3(t(D) -2)/2 \leq Volume(S^3 - K) < v_3(16t(D)
-16)$, where $v_3(\approx 1.01494)$ is the volume of a regular
hyperbolic ideal $3$-simplex.
\end{thm}

A particularly nice class of alternating diagrams is given by
$2$-bridge knots that are not torus knots.  The following lemma is a
well known consequence of work of Hatcher and Thurston.

\begin{lem} \label{2bridge}
There are $2$-bridge knots whose complements support complete
hyperbolic structures of arbitrarily large volume.
\end{lem}

\begin{pf}
It follows from \cite{HT} that $2$-bridge knots are simple and from
\cite{T} that the complement of a $2$-bridge knot that is not a torus
knot supports a complete finite volume hyperbolic structure.  

\vspace{1 mm}
\noindent
Claim: There are $2$-bridge knots that are not torus knots with
diagrams of arbitrarily high twist number.
\vspace{1 mm}

A $2$-bridge knot is determined by a sequence of integers $[c_1, \dots,
c_n]$ denoting the number of crossings in its twists, read from top to
bottom, and with $c_i$ being the number of positive crossings if $i$
is odd and the number of negative crossings if $i$ is even.  Such a
sequence gives rise to a rational number
\[\frac{p}{q} = 1 + \frac{1}{c_2 + \frac{1}{c_3 + \dots}}.\]
For instance, the Figure $8$ Knot, with sequence $[2,2]$ corresponds
to $\frac{5}{2} = 2 + \frac{1}{2}$.

Two $2$-bridge knots, with corresponding rational numbers
$\frac{p}{q}$ and $\frac{p'}{q'}$, are equivalent if and only if $p =
p'$ and $q-q'$ is divisible by $p$.  It follows from \cite{Schubert}
(for a shorter proof see \cite{Sc}) that the bridge number of a $(p,
q)$-torus knot is $min(p, q)$.  Thus a $2$-bridge knot that is also a
torus knot must be a $(2, n)$-torus knot.  The rational number
corresponding to the $(2, n)$-torus knot is $n$, an integer.  Examples
of $2$-bridge knots that are not torus knots with diagrams of
arbitrarily high twist number can thus be easily constructed, e.g.:
$[2,2]$, $[2,2,2]$, $[2,2,2,2], \dots$ These have corresponding
nonintegral rational numbers $\frac{5}{2} = 2 + \frac{1}{2}$,
$\frac{12}{5} = 2 + \frac{2}{5}$, $\frac{29}{12} = 2 + \frac{5}{12},
\dots$.

\vspace{1 mm}

Since there are $2$-bridge knots that are not torus knots with
diagrams of arbitrarily high twist number, Lackenby's Theorem \cite{L}
implies that there are $2$-bridge knots of arbitrarily large volume.
\end{pf}

\begin{defn}
A \underline{tunnel system} for a knot $K$ is a collection of disjoint
arcs $\cal T$ $= t_1 \cup \dots \cup t_n$, properly embedded in $C(K)
= S^3 - \eta(K)$ such that $C(K) - \eta(\cal T)$ is a handlebody.  The
\underline{tunnel number of K}, denoted by $t(K)$, is the least number
of arcs required in a tunnel system for $K$.

A \underline{Heegaard splitting} of a closed $3$-manifold $M$ is a
decomposition $M = V \cup_S W$ in which $V, W$ are handlebodies with
$\partial V = \partial W$ is the surface $S$, called the splitting
surface.  The \underline{genus of $M$} is the minimal genus required
for a splitting surface of $M$.
\end{defn}

The following Lemma is well known (see for instance \cite{M}).

\begin{lem} \label{tunnel1}
$2$-bridge knots have tunnel number 1.
\end{lem}

Recall the $2\pi$-Theorem (for a proof, see for instance \cite[Theorem
9]{BH}):
(Here $X(s_1, \dots, s_n)$ is the $3$-manifold obtained by Dehn filling
$X$ along $s_1 \cup \cdots \cup s_n$.)

\begin{thm} 
(Gromov-Thurston) Let $X$ be a compact orientable hyperbolic
$3$-man-\newline ifold.  Let $s_1, \dots, s_n$ be a collection of
slopes on distinct components $T_1, \dots, T_n$ of $\partial X$.
Suppose that there is a horoball neighborhood of $T_1 \cup \dots \cup
T_n$ on which each $s_i$ has length greather than $2\pi$.  Then
$X(s_1, \dots, s_n)$ has a complete finite volume Riemannian metric
with all sectional curvatures negative.
\end{thm}

More recently, in their investigation of Dehn surgery, D. Cooper and
M. Lackenby established the following relationship between the
Thurston norm of a compact hyperbolic $3$-manifold and that of its
Dehn fillings (\cite[Proposition 3.3]{CL}):

\begin{thm}
(Cooper-Lackenby) There is a non-increasing function $\beta: (2\pi,
\infty) \rightarrow (1, \infty)$, which has the following property.
Let $X$ be a compact hyperbolic $3$-manifold and let $s_1, \dots, s_n$
be slopes on distinct components $T_1, \dots, T_n$ of $\partial X$.
Suppose that there is a maximal horoball neighborhood of $T_1 \cup
\dots \cup T_n$ on which $l(s_i) > 2\pi$ for each $i$.  Then \[ \vline
X(s_1, \dots, s_n) \vline \leq \vline X \vline < \vline X(s_1, \dots,
s_n) \vline \beta(min_{1 \geq i \geq n} l(s_i)) \].
\end{thm}

Recall that the that the volume and the Thurston norm of a compact
hyperbolic $3$-manifold $M$ satisfy $vol(M) = v_3 |M|$, where $v_3$ is
the volume of a regular ideal $3$-simplex in hyperbolic $3$-space.

\begin{thm} \label{main}
There exist genus 2 closed hyperbolic 3-manifolds of arbitrarily large
volume.
\end{thm}

\begin{pf}
Let $N \epsilon {\bf R}^+$.  Choose $\varepsilon > 0$ and choose a
$2$-bridge knot $K_{p/q}$ that is not a torus knot such that its
complement $X = S^3 - \eta(K)$ has $|X| = \frac{vol(X)}{v_3} \geq
\beta(2\pi + \varepsilon)N$, for $\beta$ as provided by
Cooper-Lackenby's Theorem.  Let $r$ be a slope satisfying the
hypotheses of Thurston's Hyperbolic Surgery Theorem (see \cite[Theorem
5.8.2]{Th} or \cite[Section E.5]{BP}), then $X(r)$ is hyperbolic.

Let $\alpha$ be an arc in $X$ that is a tunnel system for $K_{p/q}$.
Let $\tilde V = \eta(\partial X \cup \alpha)$ and let $W = closure(X -
V)$.  By abusing notation slightly, we may consider $W$ to be lying in
$X(r)$.  Set $V = closure(X(r) - W)$ and $S = V \cap W$.  Then $X(r) =
V \cup_S W$ is a genus $2$ Heegaard splitting of $X(r)$.

Suppose that $V \cup_S W$ is reducible.  Then either $X(r)$ is
reducible, or $V \cup_S W$ is stabilized.  In case of the former,
$X(r)$ would be the connected sum of two lens spaces.  In case of the
latter, $X(r)$ would have genus $1$, i.e., be a Lens space, or genus $0$,
i.e., be $S^3$, but all of these outcomes would contradict the fact
that $X(r)$ is hyperbolic.  Thus $X(r)$ has genus $2$.

By the theorem of Cooper-Lackenby, \[vol(X(r)) = v_3 |X(r)| >
\frac{1}{\beta(l(r))}|X| \geq \frac{\beta(2\pi +
\varepsilon)}{\beta(l(r))}N \geq N.\]
\end{pf}

\begin{cor}
There are closed manifolds with fundamental group of rank $\leq 2$ of
arbitrarily large hyperbolic volume.
\end{cor}

\end{document}